\documentclass[10pt]{article}

\usepackage[centertags]{amsmath}
\usepackage{amsthm}
\usepackage{amssymb}
\usepackage{amscd}
\usepackage{eucal}

\newcommand{\R}{\mathbb{R}}
\newcommand{\Z}{\mathbb{Z}}

\newcommand{\iso}{\cong}           
\newcommand{\smooth}{C^\infty}
\newcommand{\suchthat}{\; : \;}

\newcommand{\id}{\mathrm{id}}

\renewcommand{\o}{\omega}

\newtheoremstyle{my}{1.5em}{0.5em}{\em}{}{\sc}{.}{0.5em}{}
\theoremstyle{my}
\newtheorem{thm}{Theorem}

\newtheorem{prop}[thm]{Proposition}
\newtheorem{cor}[thm]{Corollary}

\newenvironment{romanlist}%
{

\begin{enumerate} \parsep0em
\itemsep0em \parskip0em}{\parskip-1em \end{enumerate}}

\newcommand{\GG}{\mathcal{G}}

\title{Braids and symplectic four-manifolds\\ with abelian
fundamental group}
\author{Paul Seidel}
\date{December 18, 2001}

\begin{document}

\maketitle

\begin{abstract}
We explain how a version of Floer homology can be used as an
invariant of symplectic manifolds with $b_1>0$. As a concrete
example, we look at four-manifolds produced from braids by a surgery
construction. The outcome shows that the invariant is nontrivial;
however, it is an open question whether it is stronger than the known
ones.
\end{abstract}

On a symplectic manifold with nonzero first Betti number, there is a
distinction between symplectic and Hamiltonian vector fields. This is
usually perceived as adding to the difficulty of understanding such
manifolds, and indeed it raises many questions, some of which are
still open (for instance, the flux conjecture
\cite{lalonde-mcduff-polterovich95}). But the additional complexity
also means that more symplectic invariants become available. The aim
of this note is to draw attention to one of these, which we call
``non-Hamiltonian Floer homology''. The concept is by no means new,
but it has been underrated a bit, partially because the first
computations were in cases where it reduces to a version of ordinary
homology \cite{le-ono95}.

To demonstrate its usefulness, we consider a construction that
associates to any $d$-stranded braid (to be precise, framed spherical
transitive braid) a symplectic four-manifold with fundamental group
$\Z \times \Z/d$. This is a variation of earlier constructions due to
Smith \cite{smith01,smith00}, McMullen-Taubes
\cite{mcmullen-taubes00} and Fintushel-Stern
\cite{fintushel-stern98,fintushel-stern99}. Braids can be represented
as diffeomorphisms of a punctured two-sphere. Moreover, these
representatives can be chosen to be symplectic, and then there is a
symplectic Floer homology group measuring their fixed point theory.
We will prove that the non-Hamiltonian Floer homology of the
associated four-manifold recovers this Floer homology of the braid
(more precisely, recovers its total dimension).

At present, it is an open question whether the same information could
be obtained from Gromov-Witten theory, or even from more classical
topological invariants. Non-Hamiltonian Floer homology is not
invariant under deformations of the symplectic class, which seems to
indicate that it cannot be expressed in terms of Gromov-Witten
invariants alone. This still does not quite answer the question, so
the importance of our results remains somewhat dubious; which is one
reason why this is only an announcement, containing no proofs.

{\em Acknowledgements.} Ivan Smith has generously shared his ideas
with me; his influence on this work is more considerable than the
quotations alone would suggest. I have also had many stimulating
conversations with Stefano Vidussi.

\section{Floer homology}

Let $(M,\o)$ be a closed symplectic manifold. To any symplectic
automorphism $\phi$ of $M$ one can associate its Floer homology group
$HF_*(\phi)$, which is a finite-dimensional $\Z/2$-graded vector
space over the ``universal Novikov field'' $\Lambda$. We will use
Floer homology only for manifolds of dimension $\leq 4$, where the
definition becomes considerably easier (as a consequence of ``weak
monotonicity''). Here are some basic properties:
\begin{romanlist}
\item \label{item:arnold}
$HF_*(\id) \iso H_*(M;\Lambda)$ canonically,

\item \label{item:poincare}
$HF_*(\phi^{-1})$ is the vector space dual of $HF_*(\phi)$,

\item \label{item:conj}
for any $\phi,\psi$ there is a canonical isomorphism $HF_*(\phi) \iso
HF_*(\psi\phi\psi^{-1})$,

\item \label{item:iso}
$HF_*(\phi)$ is invariant under Hamiltonian isotopies of $\phi$.
\end{romanlist}
We emphasize that $HF_*(\phi)$ is not in general invariant under
symplectic isotopies (an exception is the case when $\id-\phi^*:
H^1(M;\R) \rightarrow H^1(M;\R)$ is an isomorphism, because then any
symplectic isotopy of $\phi$ can be replaced by a Hamiltonian isotopy
followed by conjugation).

Floer homology can be formally described as Morse theory for the
action functional on the twisted free loop space ${\mathcal
L}(M,\phi) = \{u \in \smooth(\R,M) \suchthat u(t) = \phi(u(t+1))\}$.
More concretely, it is the homology of a chain complex whose
generators are the fixed points of $\phi$ (in the generic situation
where these are nondegenerate). The grading is given by the Lefschetz
index $sign(det(1-D\phi))$. This implies that the Euler
characteristic of Floer cohomology is the Lefschetz fixed point
number:
\begin{equation} \label{eq:index}
dim\, HF_0(\phi) - dim\, HF_1(\phi) = L(\phi).
\end{equation}
For more information about the construction of Floer homology, see
\cite{dostoglou-salamon94}, \cite{seidel97}. There is a considerable
amount of additional structure on these groups, which we will not
mention at all here.

{\em An example: Floer homology for braids.} On the two-sphere $S^2$,
choose a set $\Delta = \{z_1,\dots,z_d\}$ of $d \geq 2$ marked
points. Around each of these points choose distinguished local
coordinates, which means an oriented embedding $\iota : D_\epsilon
\times \{1,\dots,d\} \rightarrow S^2$ with $\iota(0,k) = z_k$, where
$D_\epsilon$ is the closed disc of radius $\epsilon>0$. Let $\GG_d$
be the group of diffeomorphisms $\phi: S^2 \rightarrow S^2$ which
preserve $\Delta$ and are compatible with the local coordinates, in
the sense that $\phi(\iota(z,k)) = \iota(z,\sigma(k))$ for some
permutation $\sigma \in S_d$. The framed spherical braid group is
defined to be $\pi_0(\GG_d)$. An element $\beta$ of this group will
be called a transitive braid if the induced permutation of $\Delta$
is $\sigma = (1 2 \dots d)$.

To associate a Floer group to a transitive braid, choose a symplectic
form on $S^2$ such that $\iota$ is symplectic with respect to the
standard form on $D_\epsilon$. By Moser's lemma on volume forms,
$\beta$ has a symplectic representative $\phi$. Consider $S' = S^2
\setminus \iota(int(D_{\epsilon/2}) \times \{1,\dots,d\})$ and its
symplectic automorphism $\phi' = \phi\,|\,S'$. One defines
\begin{equation} \label{eq:braid}
HF_*(\beta) \stackrel{def}{=} HF_*(\phi').
\end{equation}
We are overstepping the bounds slightly, since we did not originally
introduce Floer homology for manifolds with boundary. However, this
is not really a problem: by transitivity of $\beta$, there are no
fixed points on the boundary, and an easy maximum principle argument
shows that the ``connecting orbits'' of Floer theory never touch the
boundary. (An alternative way of defining $HF_*(\beta)$ is to glue in
a torus to each boundary component, which yields a closed genus $d$
surface. One can extend $\phi'$ to a symplectic map of this surface
which permutes the tori transitively, and then define $HF_*(\beta)$
to be the Floer homology of this extension. The techniques of
\cite{seidel96b} ensure that both approaches yield the same result.)
Again using Moser's lemma, we know that the representative $\phi$ is
unique up to symplectic isotopies within $\GG_d$. All such isotopies
are Hamiltonian, even though the Hamiltonian functions involved do
not necessarily vanish near $\partial S'$. As a consequence
\eqref{eq:braid} is independent of the choice of $\phi$. This is
evident if one defines $HF_*(\beta)$ in terms of the closed genus $d$
surface, since the induced isotopy on that surface is Hamiltonian. If
one wants to stay on $S'$, the proof requires another application of
the maximum principle to solutions of the ``continuation'' equation.

It is not entirely clear how much information $HF_*(\beta)$ contains.
Conjugate braids have the same Floer homology, by \ref{item:conj}
above. A ``neck-stretching'' argument in the spirit of
\cite{seidel96b}, using the transitivity of $\beta$, shows that Floer
homology does not really see the framings: if $\tau \in \GG_d$ is a
Dehn twist along a small loop encircling just one point of $\Delta$,
then $HF_*(\beta) \iso HF_*(\beta \circ [\tau])$ for all $\beta$. On
the positive side, Nielsen theory provides a nontrivial lower bound.
Recall that for each connected component $l \in \pi_0({\mathcal
L}(S',\phi'))$ of the twisted free loop space there is a Nielsen
number $N_l(\phi') \in \Z$, which counts (with the usual sign) the
fixed points whose associated constant paths lie in that component.
Floer homology admits a corresponding splitting into direct summands,
whose Euler characteristics are the $N_l(\phi')$; this is a
refinement of \eqref{eq:index}. As a consequence
\begin{equation} \label{eq:nielsen}
dim\,HF_*(\beta) \geq \sum_l |N_l(\phi')|.
\end{equation}
This is a general feature of Floer homology, not at all limited to
braids, but it is particularly relevant in this case due to the
richness of $\pi_1(S')$. For all I know, \eqref{eq:nielsen} might be
an equality for all transitive $\beta$ (note however, that there are
plenty of counterexamples among more general surface diffeomorphisms,
starting with Poincar{\'e}'s geometric theorem).

{\em Non-Hamiltonian Floer homology.} Let $M$ be a closed symplectic
manifold. To any $a \in H^1(M;\R)$ one can associate a Hamiltonian
isotopy class of automorphisms: it consists of maps $\phi_a$ which
can be obtained from the identity by a symplectic isotopy with Calabi
class $a$ (see \cite[Chapter 10]{mcduff-salamon} for a thorough
discussion; what we call the Calabi class is the ``flux'' in their
terminology, and moreover our sign convention is opposite to theirs).
The ``non-Hamiltonian Floer homology'' $HF_*(\phi_a)$ is an invariant
of $(M,a)$. From the properties stated above one sees that for $a =
0$ it is ordinary homology; that the groups associated to $a$ and
$-a$ are dual; and that $HF_*(\phi_a) \iso HF_*(\phi_{a+\gamma})$ for
all $\gamma$ in the flux subgroup $\Gamma \subset H^1(M;\R)$ (again,
this is defined in \cite{mcduff-salamon}).

For monotone symplectic manifolds (of arbitrary dimension)
$HF_*(\phi_a)$ can be determined completely. This has been done by
L\^e and Ono in \cite{le-ono95}; their result, stated in a slightly
different form, is

\begin{thm}[L\^e-Ono] \label{th:le-ono}
Assume that $[\o] = \lambda\, c_1(M)$ with some $\lambda \neq 0$ (if
$\lambda$ is negative, they also assume that the minimal Chern number
$N$ satisfies $2N \geq dim\, M - 4$). Then
$
HF_*(\phi_a) \iso H_*(M;\underline{\Lambda}_a),
$
where $\underline{\Lambda}_a \rightarrow M$ is the flat bundle of
invertible $\Lambda$-modules canonically associated to $a$.
\end{thm}

The assumption on $N$ can probably be dropped, in view of the
technical developments which have occurred in the meantime. On the
other hand, the monotonicity condition is essential, as we will see
now.

Given an arbitrary symplectic automorphism $\psi$ of $M$, one can
form the symplectic mapping torus $E$, which is the quotient of $S^1
\times \R \times M$ by the $\Z$-action $(s,t,x) \mapsto
(s,t-1,\psi(x))$, with the obvious product symplectic structure.
Informally speaking, the purpose of mapping tori is to make maps
isotopic to the identity. In our context, this means that the
symplectic automorphism $\phi$ of $E$ given by $(s,t,x) \mapsto
(s,t,\psi(x))$ is the time one map of the symplectic vector field
$\partial/\partial t$. In terms of the notation introduced above,
$\phi = \phi_a$ with $a = [ds] \in H^1(M;\R)$.

\begin{prop} \label{th:suspension}
$HF_*(\phi) \iso HF_*(\psi) \otimes_\Lambda H_*(T^2;\Lambda)$.
\end{prop}

On a naive level, this ``suspension'' isomorphism reflects the fact
that each fixed point of $\psi$ gives rise to a $T^2$ of fixed points
of $\phi$. Note that mapping tori are never monotone: and indeed, in
contrast to Theorem \ref{th:le-ono}, non-Hamiltonian Floer homology
does not reduce to a form of ordinary homology. Another obvious
remark is that non-Hamiltonian Floer homology depends not only on
$a$, but also on the cohomology class of the symplectic form. Indeed,
if one rescales $a$ by some small amount and keeps the symplectic
form on $E$, then all fixed points of $\phi_a$ disappear; and the
same happens if one keeps $a$ fixed and changes the symplectic form
to $\delta(ds \wedge dt) + \o$, for $\delta \notin \Z$.

For the purposes of four-dimensional symplectic topology, Proposition
\ref{th:suspension} is not very satisfactory, since the mapping tori
of surface diffeomorphisms have highly nontrivial fundamental groups,
so that one can often distinguish between them on topological
grounds.

{\em Path components.} When explaining \eqref{eq:nielsen}, we had
already mentioned the decomposition of $HF_*(\phi)$ into direct
summands corresponding to components of the twisted free loop space.
For $\phi = \phi_a$ which is homotopic to the identity, this is just
the ordinary free loop space, so summands are enumerated by conjugacy
classes in $\pi_1(M)$. For our purpose, it is sufficient to have a
coarser splitting, which distinguishes only homology classes. We
denote this by
\begin{equation} \label{eq:split}
HF_*(\phi_a) \iso \bigoplus_{c \in H_1(M;\Z)} HF_*(\phi_a;c).
\end{equation}
The duality between $HF_*(\phi_a)$ and $HF_*(\phi_{-a})$ relates the
components belonging to $c$ and $-c$. Some examples: in Theorem
\ref{th:le-ono} only the summand corresponding to $c = 0$ is nonzero.
In fact, the proof of that result is essentially by deformation to $a
= 0$, where one can arrange that all fixed points lie in the trivial
connected component. In Proposition \ref{th:suspension} there can be
several nontrivial summands, but all of them are for $c \neq 0$. The
reason is that the isotopy from the identity to $\phi(s,t,x) =
(s,t,\psi(x))$ winds once around the base $T^2$ in $t$-direction.

\section{A surgery construction}

Let $\beta$ be a transitive framed spherical braid on $d$ strands. We
choose a symplectic representative $\phi$ as in the previous section.
Let $M_1 = (S^1 \times \R \times S^2)/\Z$ be the symplectic mapping
torus of $\phi$. The ``graph of the braid'' is a canonical symplectic
torus $H_1 \subset M_1$: it is the image of the embedding $T^2
\rightarrow M_1$ which sends $(s,t)$, for $k-1/d \leq t \leq k/d$, to
$(s,d\cdot t - k+1,z_k)$. Because we are working with framed braids,
the normal bundle of $H_1$ has a canonical trivialization. Moreover,
$vol(H_1) = d$. Next take $M_2 = T^4$ with the disjointly embedded
tori
\begin{align*}
 & H_1' = S^1 \times S^1 \times \{\xi_1\} \times \{\xi_2\}, \\
 & H_3' = \{\xi_3\} \times S^1 \times S^1 \times \{\xi_4\}, \\
 & H_4' = \{\xi_5\} \times S^1 \times \{\xi_6\} \times S^1
\end{align*}
where the $\xi_r$ are all different. We identify each of these tori
with $T^2$ by using coordinates in the given order. Their normal
bundles have preferred (translation-invariant) trivializations. Equip
$M_2$ with some constant symplectic form which makes all the $H_k'$
symplectic, and such that $vol(H_1') = d$. The remaining parts
$M_3,M_4$ will be elliptically fibered K3 surfaces, with embedded
tori $H_3 \subset M_3$, $H_4 \subset M_4$ which are the fibres. We
identify them with $T^2$ in an arbitrary way, and use the standard
trivialization of their normal bundles (given by the fibration). The
symplectic forms should be normalized in such a way that $vol(H_3) =
vol(H_3')$, $vol(H_4) = vol(H_4')$. Now glue together all these
pieces by Gompf-style sums, pairing the tori $H_k$ with $H_k'$ for $k
= 1,3,4$. The outcome (keeping the choices of $M_2,M_3,M_4$ fixed)
depends up to symplectic isomorphism only on the conjugacy class of
the framed braid $\beta$. We denote it by $M^\beta$.

The idea can be summarized as follows. By deforming $\phi$ to the
identity inside $Diff(S^2)$, one can identify $M_1$ with $T^2 \times
S^2$. From this point of view, the two-torus $H_1 \subset T^2 \times
S^2$ is knotted in a way which is determined by the braid (as
observed in \cite{fintushel-stern99} and \cite{smith01}, one can use
the fundamental group of the complement to verify that many different
knot types occur). At this point, one could choose to directly glue
in a K3 surface to $H_1$. This is an appealing possibility, somewhat
similar to \cite[Section 5]{fintushel-stern99}, but since the
resulting fundamental group is $\Z/d$, non-Hamiltonian Floer homology
cannot be applied in a meaningful way. The role of the intermediate
piece $M_2$, which we have borrowed from \cite{mcmullen-taubes00} and
\cite{smith00}, is to let a slightly larger part of $\pi_1(M_1
\setminus H_1)$ survive.

{\em Topological aspects.} As has been already mentioned,
$\pi_1(M^\beta) \iso \Z \times \Z/d$ for all $\beta$. More
explicitly, let $z_0 \in S^2 \setminus \Delta$ be a fixed point of
$\phi$, and consider the loops $l_1 = S^1 \times \{0\} \times
\{z_0\}$, $l_2 = \{0\} \times S^1 \times \{z_0\}$ in $M_1 \setminus
H_1 \subset M^\beta$, oriented in the obvious way. The isomorphism
can be chosen such that $[l_1], [l_2]$ correspond to $(1,0)$ and
$(0,1)$ in $\Z \oplus \Z/d$, respectively. To verify this, it is
useful to carry out the gluing in a particular order. Consider the
manifold $M_2'$ obtained by putting together $M_2,M_3,M_4$ in the way
described above; this still contains the torus $H_1'$. It is a
familiar fact that $\pi_1(M_3 \setminus H_3) = \pi_1(M_4 \setminus
H_4) = 1$, and as a consequence
\begin{align*}
 \pi_1(M_2' \setminus H_1')
 & \iso \pi_1(M_2 \setminus H_1') /
 \langle \pi_1(H_3'), \pi_1(H_4') \rangle \\
 & \iso \langle d_1,d_2,d_3,d_4 \suchthat [d_1,d_2] = 1 \rangle /
 \langle d_2,d_3,d_4 \rangle \\
 & = \langle d_1 \rangle \iso \Z.
\end{align*}
Note that $d_1$ is the first longitude of $H_1'$; the other longitude
(which would be $d_2$) and the meridian (which would be $[d_3,d_4])$
have got killed. As for the remaining piece, the fundamental group of
$M_1 \setminus H_1$ is quite large, but it is generated by $l_1$
(which is the first longitude of $H_1$, and commutes with all other
elements) together with $l_2$ and various conjugates of the meridian
of $H_1$. Joining together this with $M_2' \setminus H_1'$ kills the
meridian and identifies $l_1$ with $d_1$, from which one sees that
the fundamental group becomes abelian.

The characteristic numbers are $c_2(M^\beta) = 48$, $c_1(M^\beta)^2 =
0$. In fact $-d \cdot c_1(M^\beta)$ can be represented by the
disjoint union of $6d-2$ embedded symplectic tori, each of which has
trivial normal bundle: $2d-2$ parallel copies of $H_1$, and $2d$
copies of $H_3$ and $H_4$ each (see \cite{smith00} for how to do this
kind of computation). The next step would be to compute the homology
of the universal cover, as a $\pi_1(M^\beta)$-module, and the
intersection form on it (presumably, that goes a long way towards
determining the homeomorphism type of $M^\beta$). We have not done
this, but informal considerations suggest that it might turn out to
be the same for all $\beta$.

{\em Floer homology.} Let $a_1 \in H^1(M^\beta;\R)$ be the unique
class with $\langle a_1,[l_1] \rangle = 1$. The next result
determines certain summands in the splitting \eqref{eq:split} of the
non-Hamiltonian Floer group $HF_*(\phi_{a_1})$.

\begin{thm} \label{th:main}
$HF_*(\phi_{a_1};[l_2]) \iso HF_*(\beta) \otimes_\Lambda
H_*(T^2;\Lambda)$. Moreover $HF_*(\phi_{a_1};c) = 0$ for any other
nonzero torsion class $c \in H_1(M^\beta;\Z)$.
\end{thm}

Consider the direct sum of the groups $HF_*(\phi_a;c)$ where $a$
ranges over the two generators of $H^1(M^\beta;\Z) \subset
H^1(M^\beta;\R)$, and $c$ over all nonzero torsion elements of
$H_1(M^\beta;\Z)$. Using the duality between $HF_*(\phi_a;c)$ and
$HF_*(\phi_{-a};-c)$ one computes that the total dimension of the
direct sum is $8 \,dim( HF_*(\beta))$. On the other hand, the direct
sum is defined without reference to any particular basis of homology.
As a consequence:

\begin{cor}
The total dimension of $HF_*(\beta)$ is a symplectic invariant of
$M^{\beta}$.
\end{cor}

As in our discussion of mapping tori, $HF_*(\phi_a;c)$ becomes zero
if one changes the symplectic class slightly, by changing the area of
the $T^2$ factor in $M_1$, and rescaling the symplectic forms of the
other $M_k$ accordingly. However, one can compensate for this by
rescaling $a$, and recover the Floer homology groups in this way. Of
course, it is not clear what happens under ``large'' deformations of
the symplectic class.

The proof of Theorem \ref{th:main} consists of two steps. The first
is a variant of Proposition \ref{th:suspension} adapted to surfaces
with boundary. The second step is a ``Mayer-Vietoris'' argument in
which one considers the behaviour of Floer homology under Gompf sums.
In general this is a hard problem, as a look at the formulae for
Gromov-Witten invariants shows \cite{ionel-parker98, li-ruan98}, but
the particular case needed here is fairly simple. The reason is
essentially topological: the components $c$ which we are interested
in lie in the image of $H_1(M_1 \setminus H_1;\Z) \rightarrow
H_1(M^\beta;\Z)$, but not in that of $H_1(M_2' \setminus H_1';\Z)
\rightarrow H_1(M^\beta;\Z)$. Therefore only the fixed points lying
in $M_1 \setminus H_1$ are relevant, which are precisely those coming
from the braid.

Finally, we would like to point out that non-Hamiltonian Floer
homology also merits some attention in higher dimensions. For
instance, take a symplectic manifold $M$ with an automorphism $\phi$
that is differentiably, but not symplectically, isotopic to the
identity (in dimension $\geq 4$, plenty of such exist). Then the
symplectic mapping torus $E$ is diffeomorphic to $T^2 \times M$ but
with a potentially nonstandard symplectic structure, which one could
try to detect using Proposition \ref{th:suspension}. There are also
examples of ``fragile'' symplectic automorphisms, which become
symplectically isotopy to the identity after a slight change of the
symplectic class \cite{seidel97}. In that case, one can hope to show
that $E$ is symplectically deformation equivalent, but not
symplectomorphic, to $T^2 \times M$.

\bibliographystyle{amsplain}

\begin{thebibliography}{1}

\bibitem{dostoglou-salamon94}
S.~Dostoglou and D.~Salamon, \emph{Self dual instantons and
holomorphic
  curves}, Annals of Math. \textbf{139} (1994), 581--640.

\bibitem{fintushel-stern98}
R.~Fintushel and R.~Stern, \emph{Knots, links, and $4$-manifolds},
Invent.
  Math. \textbf{134} (1998), 363--400.

\bibitem{fintushel-stern99}
\bysame, \emph{Symplectic surfaces in a fixed homology class}, J.
Differential Geom.
  \textbf{52} (1999), 203--222.

\bibitem{ionel-parker98}
E.-N. Ionel and T.~Parker, \emph{Gromov-{W}itten invariants of
symplectic
  sums}, Math. Res. Lett. \textbf{5} (1998), 563--576.

\bibitem{lalonde-mcduff-polterovich95}
F.~Lalonde, D.~McDuff, and L.~Polterovich, \emph{On the flux
conjectures},
  Geometry, topology, and dynamics (Montreal, PQ, 1995) (F.~Lalonde, ed.),
  Amer. Math. Soc., 1998, pp.~69--85.

\bibitem{le-ono95}
{Le Hong Van} and K.~Ono, \emph{Cup-length estimate for symplectic
fixed
  points}, Contact and symplectic geometry (C.~B. Thomas, ed.), Cambridge Univ.
  Press, 1996, pp.~268--295.

\bibitem{li-ruan98}
A.-M. Li and Y.~Ruan, \emph{Symplectic surgery and {G}romov-{W}itten
invariants
  of {C}alabi-{Y}au 3-folds}, Invent. Math. \textbf{145} (2001), 151--218.

\bibitem{mcduff-salamon}
D. McDuff and D.~Salamon, \emph{{I}ntroduction to symplectic
topology}, Oxford University Press, 1995.

\bibitem{mcmullen-taubes00}
C. McMullen and C.~Taubes, \emph{{4}-manifolds with inequivalent
symplectic forms and 3-manifolds with inequivalent fibrations}, Math.
Res. Lett. \textbf{6} (1999), 681--696.

\bibitem{seidel96b}
P.~Seidel, \emph{The symplectic {F}loer homology of a {D}ehn twist},
Math.
  Research Lett. \textbf{3} (1996), 829--834.

\bibitem{seidel97}
\bysame, \emph{Floer homology and the symplectic isotopy problem},
Ph.D.
  thesis, Oxford University, 1997.

\bibitem{smith00}
I.~Smith, \emph{On moduli spaces of symplectic forms}, Math. Res.
Lett.
  \textbf{7} (2000), 779--788.

\bibitem{smith01}
\bysame, \emph{Symplectic submanifolds from surface fibrations},
Pacific J. Math. \textbf{198} (2001), 197--205.

\end{thebibliography}
\providecommand{\bysame}{\leavevmode\hbox
to3em{\hrulefill}\thinspace}

Paul Seidel \\
Centre de Math\'ematiques \\
U.M.R. 7640 du C.N.R.S. \\
Ecole Polytechnique \\
91128 Palaiseau cedex \\
France

\end{document}